\documentclass[10pt,leqno]{article}
\usepackage{amssymb,amsbsy,amsmath,amsfonts,amssymb,amscd, mathrsfs}

\usepackage[english]{babel}
\usepackage[T1]{fontenc}
\usepackage{indentfirst}

\makeatletter
\@addtoreset{equation}{section}
\makeatother

\newtheorem{statement}{}[section]
\newtheorem{theoreme}[statement]{Theorem}
\newtheorem{lemme}[statement]{Lemma}

\newtheorem{proposition}[statement]{Proposition}

\newtheorem{corollaire}[statement]{Corollary}

\newcommand\C{\mathbb C}

\newcommand\R{\mathbb R}
\newcommand\T{\mathbb T}
\newcommand\D{\mathbb D}
\newcommand\Z{\mathbb Z}
\newcommand\e{{\rm e}}

\newcommand\eps{\varepsilon}
\newcommand\ind{{\rm 1\kern-.30em I}}
\newcommand\qed{\hfill $\square$}
\renewcommand \Re{{\mathfrak R}{\rm e}\,}
\renewcommand \Im{{\mathfrak I}{\rm m}\,}

\title{\bf Compact composition operators on $H^2$ and Hardy-Orlicz spaces}
\author{\it Pascal Lef\`evre, Daniel Li,\\ \it Herv\'e Queff\'elec, Luis Rodr{\'\i}guez-Piazza}

\date{\footnotesize \today}

\begin{document}

\maketitle

\noindent{\bf Abstract.} \emph{We compare the compactness of composition 
operators on $H^2$ and on Orlicz-Hardy spaces $H^\Psi$. We show in particular that exists an Orlicz function 
$\Psi$ such that $H^{3+\eps} \subseteq H^\Psi \subseteq H^3$ for every $\eps >0$, and a composition 
operator $C_\phi$ which is compact on $H^3$ and on $H^{3+\eps}$, but not compact on $H^\Psi$.}
\medskip

\medskip



\section{Introduction}

Composition operators on the Hardy space $H^2$ have given rise to a lot of papers since the beginning of 
the seventies (see \cite{Shap-livre} and \cite{Co-McC}, and references therein, for an overview until the middle 
of the eighties). In particular, criteria of compactness (\cite{Shap}, \cite{McCluer}) or of membership in  
Schatten classes (\cite{Lue} \cite{Lue-Zhu}) were found. But only few concrete examples of such operators 
were known (see \cite{Shap-Tay}, \cite{Carroll-Cowen}, \cite{Jones}, \cite{Zhu}). In \cite{JFA}, we constructed 
explicit examples of compact composition operators $\C_\phi \colon H^2 \to H^2$ which are, or not, in some 
Schatten classes $S_p$. On the other hand, we study in \cite{CompOrli} composition operators on Hardy-Orlicz 
spaces $H^\Psi$, and characterize their compactness. In this paper, we shall continue our investigation of 
composition operators on Hardy-Orlicz spaces, and study for which Orlicz functions $\Psi$ the explicit 
examples of \cite{JFA} are compact.
\medskip

Let us describe more precisely the content of this paper.\par
In Section~\ref{Delta indice deux} and Section~\ref{Delta super deux}, we consider the Hardy-Orlicz 
space $H^\Psi$ besides the Hardy space $H^2$, and we prove that for every Orlicz function $\Psi$ which does 
not satisfy the condition $\Delta_2$, there exists a symbol $\phi$ such that $C_\phi$ is compact on $H^2$, 
but not compact on $H^\Psi$ (Theorem~\ref{non psi-compact}). When the Orlicz function $\Psi$ grows 
fast enough, we can ensure, beyond the non compactness of $C_\phi$ on $H^\Psi$, that, for every $p > 2$, 
$C_\phi \colon H^2 \to H^2$ is moreover in the Schatten 
class $S_p$ (Proposition~\ref{croissance Psi}). But if $\Psi$ does not grow too fast, this is not possible: if 
$C_\phi \in S_p$ for some $p > 0$, we must have the compactness of $C_\phi$ on $H^\Psi$ 
(Proposition~\ref{Proposition Psi pas tres croissante}). When $\Psi$ grows very fast, \emph{i.e.} satisfies the 
condition $\Delta^2$, we show that $C_\phi$ may be in $S_p$ for every $p > 0$, although $C_\phi$ is 
not compact on $H^\Psi$ (Proposition~\ref{proposition Schatten et Delta super deux}); conversely, if 
$C_\phi$ is compact on $H^\Psi$, then $C_\phi$ must be in $S_p$ for every $p > 0$ 
(Proposition~\ref{Delta super deux et Carleson}). \par
In Section~\ref{phi indice theta dans H Psi}, we then characterize (Theorem~\ref{Psi delta 1}) the Orlicz 
functions $\Psi$ for which the composition operator $C_{\phi_\theta}$, which appeared in 
\cite{JFA}, Section~5, is compact on $H^\Psi$. When $\Psi$ grows fast enough (in  
particular if $\Psi$ satisfies the so-called condition $\Delta^1$), then $C_{\phi_\theta}$ is not compact on 
$H^\Psi$ (Theorem~\ref{meilleure condition sur Psi}); this gives, in particular, an improvement of 
Proposition~\ref{croissance Psi}. This characterization allows us to get the following striking 
result (Theorem~\ref{pas d'interpolation}): there exists an Orlicz function $\Psi$ such that 
$H^{3+\eps} \subseteq H^\Psi \subseteq H^3$ for every $\eps >0$, and a composition operator $C_\phi$ 
which is compact on $H^3$ and on $H^{3+\eps}$, but not compact on $H^\Psi$.\par 
\bigskip

\section{Notation}

We shall denote by $\D$ the open unit disk of the complex plane: $\D = \{ z\in \C\,; \ |z| < 1\}$, and by 
$\T = \partial \D$ its boundary: $\T =\{ z\in \C\,;\ |z | = 1\}$. We shall denote by $m$ the normalized 
Lebesgue measure on $\T$.\par
For every analytic self-map $\phi \colon \D \to \D$, the composition operator $C_\phi$ is the map 
$f \mapsto f\circ \phi$. By Littlewood's subordination principle (see \cite{Duren}, Theorem~1.7), every 
composition operator maps continuously every Hardy space $H^p$ ($p >0$) into itself, as 
well as (\cite{CompOrli}, Proposition~3.12) the Hardy-Orlicz spaces $H^\Psi$ (see below for their 
definition).\par
\smallskip

For every $\xi \in \T$ and $0 < h < 1$, the Carleson window $W (\xi, h)$ is the set
\begin{displaymath}
W (\xi, h) = \{ z \in \D\,;\ |z| \geq 1 - h\quad \text{and} \quad |\arg ( z \bar{\xi})| \leq h\}.
\end{displaymath}
For every finite positive measure $\mu$ on $\D$, one sets:
\begin{displaymath}
\rho_\mu (h) = \sup_{\xi \in \T} \mu[W (\xi, h)].
\end{displaymath}
We shall call this function $\rho_\mu$ the \emph{Carleson function} of $\mu$.\par
When $\phi \colon \D \to \D$ is an analytic self-map of $\D$, and $\mu =m_\phi$ is the pullback 
measure defined on $\D$, for every Borel set $B \subseteq \D$, by:
\begin{displaymath}
m_\phi (B) = m (\{\xi\in \T\,;\ \phi^\ast (\xi) \in B \}),
\end{displaymath}
where $\phi^\ast$ is the boundary values function of $\phi$, we shall denote $\rho_{m_\phi}$ by 
$\rho_\phi$. In this case, we shall say that $\rho_\phi$ is the Carleson function of $\phi$.\par
For $\alpha \geq 1$, we shall say that $\mu$ is an \emph{$\alpha$-Carleson measure} if 
$\rho_\mu (h) \lesssim h^{\alpha}$. For $\alpha =1$, $\mu$ is merely said to be a \emph{Carleson measure}. 
\par
The Carleson Theorem (see \cite{Duren}, Theorem~9.3) asserts that, for $1 \leq p < \infty$ (actually, for 
$ 0 < p < \infty$), the canonical inclusion $j_\mu \colon H^p \to L^p (\mu)$ is bounded if and only if 
$\mu$ is a Carleson measure. Since every composition operator $C_\phi$ is continuous on $H^p$, it defines a 
continuous map $j_\phi \colon H^p \to L^p (\mu_\phi)$; hence every pull-back measure $\mu_\phi$ is a 
Carleson measure.\par
When $C_\phi \colon H^2 \to H^2$ is compact, one has, as it is easy to see:
\begin{equation}\label{condition phi plus petit que un au bord}
|\phi^\ast| < 1 \qquad \text{\emph{a.e.} on } \partial \D.
\end{equation}
Hence, we shall only consider in this paper 
symbols $\phi$ for which \eqref{condition phi plus petit que un au bord} is satisfied (which is the case, as we 
said, when $C_\phi$ is compact on $H^2$). Recall (\cite{CompOrli}, Theorem~4.3) that, for every Orlicz 
function $\Psi$, as defined below,  $C_\phi \colon H^2 \to H^2$ is compact whenever 
$C_\phi \colon H^\Psi \to H^\Psi$ is.\par
\smallskip

B. MacCluer (\cite{McCluer}, Theorem~1.1) has shown, assuming 
condition~\eqref{condition phi plus petit que un au bord}, that $C_\phi$ is compact on $H^p$ if and 
only if $\rho_\phi (h) = o\, (h)$, as $h$ goes to $0$.
\bigskip

Throughout this paper, the notation $f\approx g$ will mean that there are two constants $0 < c < C < +\infty$ 
such that $ c f(t) \leq g(t) \leq C f(t)$ (for $t$ sufficiently near of a specified value), and the notation 
$ f (t) \lesssim g (t)$, when $t$ is in the neighbourhood of some value $t_0$, will have the same meaning as  
$g = O (f)$.\par
\medskip

Note that, in this paper, we shall not work, most often, with exact inequalities, but with 
inequalities up to constants. It follows that we shall not actually work with true Carleson windows 
$W (\xi, h)$, but with distorted Carleson windows: 
\begin{displaymath}
\tilde W (\xi, h) =\{ z\in \D\,;\ |z| \geq 1 - ah \quad \text{and} \quad | \arg (z \bar{\xi} )|\leq bh\},
\end{displaymath}
where $a,b > 0$ are given constants. Since, for a given symbol $\phi$, one has:
\begin{displaymath}
m_\phi \big(W (\xi, c\,h) \big) \leq m_\phi \big(\tilde W (\xi, h) \big) \leq m_\phi \big(W (\xi, C\,h) \big) 
\end{displaymath}
for some constants $c= c(a,b)$ and $C = C (a,b)$ which only depend on $a$ and $b$, that will not matter 
for our purpose.
\medskip

Let us now recall the definition of the Hardy-Orlicz spaces.\par
An Orlicz function is a non-decreasing convex function $\Psi\colon [0, \infty] \to [0,\infty]$ such that 
$\Psi(0) = 0$ and $\Psi(\infty) =\infty$. To avoid pathologies, we assume that the Orlicz function $\Psi$ 
has the following additional properties: $\Psi$ is continuous at $0$, strictly convex (hence increasing), and 
such that
\begin{displaymath}
\frac{\Psi(x)}{x}\mathop{\longrightarrow}_{x\to \infty} \infty.
\end{displaymath}
This is essentially to exclude the case of $\Psi(x) = ax$.\par
The Orlicz space $L^\Psi (\T)$ is the space of all (equivalence classes of) measurable functions 
$f\colon \T\to\C$ for which there is a constant $C > 0$ such that:
\begin{displaymath}
\int_\T \Psi\bigg(\frac{\vert f (t)\vert} {C}\bigg)\,dm(t) < +\infty
\end{displaymath}
and then $\Vert f\Vert_\Psi$, the \emph{Luxemburg norm}, is the infimum of all possible
constants $C$ such that this integral is $\leq 1$. The Hardy-Orlicz space $H^\Psi$ is the space of all 
$f \in H^1$ such that the boundary values function $f^\ast$ of $f$ is in $L^\Psi (\T)$. We refer to 
\cite{CompOrli}, Section~3  for more details.\par


\section{Orlicz functions without condition $\Delta_2$} \label{Delta indice deux}

In \cite{CompOrli}, Corollary~3.26, we showed that for every Orlicz function $\Psi$ which satisfies the 
growth condition $\Delta^2$, there exists a symbol $\phi\colon\D \to \D$ which induces a compact 
composition on $H^2$, but a non-compact composition operator on $H^\Psi$. We shall generalize this 
below. Recall that an Orlicz function $\Psi$ satisfies the $\Delta_2$ condition if there exists a constant 
$C > 0$ such that $\Psi (2x) \leq C \Psi (x)$ for $x$ large enough.

\begin{theoreme}\label{non psi-compact}
For every Orlicz function $\Psi$ which does not satisfy the condition $\Delta_2$, there exists an analytic 
self-map $\phi \colon \D \to \D$ such that $C_\phi$ is compact on $H^2$ but not compact on $H^\Psi$. 
\end{theoreme}

\noindent{\bf Proof.} $\Psi \notin \Delta_2$ means that 
$\limsup\limits_{x\to \infty} \frac{\Psi (2x)}{\Psi (x)} =+\infty$. We can hence find a sequence $(x_n)$, 
which increases to $+\infty$, such that $\frac{\Psi (2x_n)}{\Psi (x_n)}$ increases to $+\infty$.\par
Set $h_n= \frac{1}{\Psi(x_n)}$ and $c_n = \frac{\Psi (x_n)}{\Psi (2x_n)}\,$. By construction, $(h_n)_n$ 
and $(c_n)_n$ decrease to $0$. We may assume that $c_n < \pi$, $n\geq 1$.\par
We are going to construct an analytic function $\phi\colon \D \to \D$ such that: 
\begin{equation}\label{condition petit o}
\rho_\phi (h) = o\, (h),
\end{equation}
but:
\begin{equation}\label{pas R0}
\hskip 2cm \rho_\phi (h_n) \geq c_n h_n\,,\qquad  \text{for all } n\geq 1.
\end{equation} 
Condition \eqref{condition petit o} will ensure that $C_\phi\colon H^2 \to H^2$ is compact, by MacCluer's 
theorem, and condition \eqref{pas R0}, which reads:
\begin{displaymath}
\hskip 2cm \rho_\phi (h_n) \geq \frac{1}{\Psi \big( 2 \Psi^{-1} (1/h_n)\big)}\,\raise 1,5pt \hbox{,}
\qquad  \text{for all } n\geq 1,
\end{displaymath}
will ensure (by \cite{CompOrli}, Theorem~4.11), that $C_\phi\colon H^\Psi \to H^\Psi$ is not compact.\par
\medskip

For that purpose, we shall use the ``general construction'' made in \cite{JFA}, \S~3.2. Let us recall this 
construction.\par
Let
\begin{equation}\label{forme de f}
f (t) = \sum_{k=0}^\infty a_k \cos (kt)
\end{equation}
be an even, non-negative, $2\pi$-periodic continuous function, vanishing at the origin: $f (0) = 0$, and such 
that:
\begin{equation}\label{deuxieme condition}
\text{$f$ is \emph{strictly} increasing on $[0,\pi]$.}
\end{equation}
The Hilbert transform (or conjugate function) ${\cal H}f$ of $f$ is:
\begin{displaymath}
{\cal H}f (t) =\sum_{k=1}^\infty a_k \sin (kt).
\end{displaymath}
We shall assume moreover that, as $t$ tends to zero:
\begin{equation}\label{derivee de Hf}
\big({\cal H}f\big)' (t) = o\, ( 1/t^2).
\end{equation}
Let now $F \colon \D \to \Pi^+ = \{ \Re z >0\}$ be the analytic function whose boundary values are:
\begin{equation}\label{def grand F}
F^\ast (\e^{it}) = f (t) + i {\cal H}f (t).
\end{equation}
One has:
\begin{displaymath}
\qquad F (z) = \sum_{k=0}^\infty a_k z^k, \qquad |z| < 1,
\end{displaymath}
and we define:
\begin{equation}\label{def grand Phi}
\qquad \Phi (z) = \exp \big( - F(z)\big), \qquad z\in \D.
\end{equation}
Since $f$ is non-negative, one has:
\begin{displaymath}
\Re F (z) = \frac{1}{2\pi} \int_{-\pi}^{\pi} f (t) P_z (t)\,dt > 0,
\end{displaymath}
$P_z$ being the Poisson kernel at $z$, 
so that $| \Phi (z) | < 1$ for every $z \in \D$: $\Phi$ is an analytic self-map of $\D$, and 
$|\Phi^\ast | = \exp (- f ) < 1$ \emph{a.e.} . Note that the condition $f (0) = 0$ means that $\Phi^\ast (1) = 1$, 
so that $\| \Phi \|_\infty =1$.\par
We then perturb $\Phi$ by considering:
\begin{equation}\label{M}
\qquad M (z) = \exp\Big( - \frac{1+ z }{ 1 -z }\Big)\,, \qquad |z| <1,
\end{equation}
and 
\begin{equation}\label{def general phi}
\qquad \phi ( z) = M (z) \Phi (z), \qquad | z | <1.
\end{equation}\label{lemme 3.6}
We proved in \cite{JFA}, Lemma~3.6, that 
\begin{equation}
h f^{-1} (h) \lesssim \rho_\phi (h) \lesssim h f^{-1} (2h)
\end{equation}
when $h$ tends to $0$.
\medskip

We shall now rely on the following simple lemma, whose proof is postponed.

\begin{lemme}\label{lemme construction fonction}
Let $(h_n)$ and $(c_n)$ be two decreasing vanishing sequences of positive numbers, with 
$0< h_n< \pi$, $0 < c_n < \pi$.\par
There exists an even $2\pi$-periodic ${\cal C}^2$ function $f \colon \R \to \R$ such that: \par
1) $f (0) = 0$ and $f$ is strictly increasing on $[0,\pi]$;\par
2) $f (c_n ) \leq h_n$, for all $n \geq 1$.
\end{lemme}

We are in the situation of the general construction, and by\eqref{lemme 3.6}, we have:
\begin{displaymath}
\rho_\phi (h) = O \big( h f^{-1} (2h) \big) = o\, (h),
\end{displaymath}
since $f^{-1}(2h) \to 0$ as $h \mathop{\to}\limits^{\scriptscriptstyle >} 0$. MacCluer's criterion implies that  
$C_\phi \colon H^2 \to H^2$ is compact. On the other hand, \eqref{lemme 3.6}  again implies that
\begin{displaymath}
\rho_\phi (h_n ) \gtrsim h_n f^{-1} (h_n) \gtrsim c_n h_n,
\end{displaymath}
by Lemma~\ref{lemme construction fonction}.\par
Therefore \eqref{pas R0} holds, and this ends the proof of Theorem~\ref{non psi-compact}. \qed
\medskip

\noindent{\bf Remark.} Actually, we do not use that $f$ is ${\cal C}^2$, but only that it is continous. However, 
if $f\in {\cal C}^2$, then $\widehat{f ''} \in \ell_2 (\Z)$ and it follows, by the Cauchy-Schwarz inequality, 
since $\widehat{f ''}(k) = (ik)^2 \hat f (k)$, that $\sum_{k=0}^\infty k\,|\hat f (k)| < +\infty$; hence  
both $f$ and ${\cal H}f$ are ${\cal C}^1$ functions, and we proved in \cite{JFA}, Lemma~3.5, 
that then the Carleson function $\rho_\Phi$ of $\Phi$ is not $o\,(h)$ when $h$ goes to $0$, and so the 
composition operator $C_\Phi \colon H^2 \to H^2$ is not compact, although $|\Phi^\ast| = |\phi^\ast|$ 
on $\partial \D$ and $C_\phi \colon H^2 \to H^2$ is compact.
\bigskip

\noindent{\bf Proof of Lemma~\ref{lemme construction fonction}.} Let $h_0 = c_0 =\pi$, and 
$\varphi \colon [0, \pi] \to \R_+$ be a step function such that:
\begin{equation}\label{step function}
\int_0^{c_n} \varphi (t)\,dt = h_n\,, \qquad n=1,2,\ldots
\end{equation}
We can take $\varphi (t) = \frac{h_j - h_{j+1}}{c_j - c_{j+1}}$ for $ t \in (c_{j+1}, c_j)$, $j= 0,1, \ldots$, since 
then:
\begin{align*}
\int_0^{c_n} \varphi (t)\,dt 
& = \sum_{j=n}^\infty \int_{c_{j+1}}^{c_j} \varphi (t)\,dt = 
 \sum_{j=n}^\infty (c_j - c_{j+1})  \frac{h_j - h_{j+1}}{c_j - c_{j+1}} \\
& = \sum_{j=n}^\infty (h_j - h_{j+1}) = h_n.
\end{align*}
Let now $f \colon \R \to \R$ be the even $2\pi$-periodic function whose values on $[0,\pi]$ are given by:
\begin{displaymath}
f (t) = \frac{1}{\pi^3} \int_0^t (t - u)^3 \varphi (u)\,du.
\end{displaymath}
This function is clearly ${\cal C}^2$, with
\begin{displaymath}
f '' (t) = \frac{6}{\pi^3} \int_0^t (t -u) \varphi (u)\,du,
\end{displaymath}
and $ f( 0) = f ' (0) = f '' (0) =0$; it is strictly increasing on $[0, \pi]$, and
\begin{displaymath}
f (c_n) = \frac{1}{\pi^3} \int_0^{c_n} (c_n - u)^3 \varphi (u)\,du 
\leq \frac{c_n^3}{\pi^3} \int_0^{c_n} \varphi (u)\,du \leq \int_0^{c_n} \varphi (u)\,du = h_n,
\end{displaymath}
due to \eqref{step function}. Therefore $f$ fulfills all the requirements of 
Lemma~\ref{lemme construction fonction}. \qed
\bigskip

\noindent{\bf Remark.} It is not clear whether one can ensure in Theorem~\ref{non psi-compact} that 
the composition operator is moreover in some Schatten class $S_p$ for $p < \infty$. It would be the case, by 
\cite{JFA}, Corollary~3.2, if one can construct a symbol $\phi$ such that its pull-back measure $m_\phi$ 
is $\alpha$-Carleson for some $\alpha > 1$. This can be obtained if we can construct, in 
Lemma~\ref{lemme construction fonction}, our function $f$ such that $f^{(N)} (0) \not= 0$ for some 
$N > 1$, since then $f (h) \approx h^N$ and $f^{-1}(h) \approx h^{1/N}$, and $\phi$ would be 
$\alpha$-Carleson with $\alpha = 1 + \frac{1}{N}$, which would imply that $C_\phi \in S_p$ for 
$p > 2N$. However, in general, we cannot obtain that $f^{(N)} (0) \not= 0$ for some $N > 1$. Indeed, 
take $\Psi (x) = \e^{(\log (x+1) )^2} - 1$; then 
$h_n = \frac{1}{\Psi (x_n)} \approx \exp \big[- (\log x_n)^2\big]$, whereas, since
\begin{displaymath}
\Psi (2x) \approx \exp\big[ (\log x)^2\big]\,\exp (2\log x),
\end{displaymath}
one has $c_n =\frac{\Psi (x_n)}{\Psi (2 x_n)} \approx \exp (- 2 \log x_n)$, so that 
$h_n = o\,(c_n^N)$ for every $N < \infty$. Since $f (c_n) \leq h_n$, all the possible derivatives of $f$ at $0$ 
must be zero.\par
\bigskip

Besides, when $\Psi$ does not grow too fast, one has the following result.

\begin{proposition}\label{Proposition Psi pas tres croissante}
Let $\Psi$ be an Orlicz function such that, for every $A > 1$, one has, as $x$ goes to infinity:
\begin{equation}\label{autre condition sur Psi}
\Psi (Ax) = O \big( \Psi (x) \big(\log \Psi (x) \big)^\eps \big)\,, \qquad \forall \eps >0.
\end{equation}
Then, if  $C_\phi \in S_p$ for some $p > 0$, the composition operator $C_\phi \colon H^\Psi \to H^\Psi$ is 
compact.
\end{proposition}

\noindent{\bf Remark.} If the condition~\eqref{autre condition sur Psi} is satisfied for some $A >1$, it is 
actually satisfied for every $A > 1$. In fact, $\Psi (Ax) = O \big( \Psi (x) \big(\log \Psi (x) \big)^\eps \big)$ 
implies that $\Psi (A^2 x) = O \big( \Psi (x) \big(\log \Psi (x) \big)^{2 \eps} \big)$, and we can replace $A$ 
by $A^2$, $A^4$, $A^8$, \dots Note also that we can replace the \emph{big-oh} assumption by a 
\emph{little-oh} one: if the proposition holds with a \emph{little-oh} condition, it also holds with a 
\emph{big-oh} one, because, for every $\eps' > 0$, if one has \eqref{autre condition sur Psi}, then, with 
$\eps < \eps'$:
\begin{displaymath}
\frac{\Psi (A x)}{\Psi (x) \big(\log \Psi (x) \big)^{\eps'}} \leq K_\eps\,\big( \log \Psi (x) \big)^{\eps - \eps'} 
\mathop{\longrightarrow}_{x \to \infty} 0.
\end{displaymath}
\medskip

\noindent{\bf Proof.} If $C_\phi \in S_p$, then, by \cite{JFA}, Proposition~3.4, one has:
\begin{displaymath}
\rho_\phi (h) = o\,\bigg( \frac{h}{[\log (1/h)]^\delta}\bigg)\,,
\end{displaymath}
with $\delta =2/p > 0$. Then, \eqref{autre condition sur Psi} with $x = \Psi^{-1} (1/h)$, gives:
\begin{displaymath}
\rho_\phi (h) = o\,\bigg( \frac{1}{\Psi (x) [\log \Psi (x) ]^\delta}\bigg) 
= o\,\bigg( \frac{1}{\Psi (Ax)}\bigg) = o\,\bigg( \frac{1}{\Psi \big( A \Psi^{-1} (1/h)\big) }\bigg)\,,
\end{displaymath}
which implies, by \cite{CompOrli}, Theorem~4.18, that $C_\phi$ is compact on $H^\Psi$. \qed
\bigskip

For example, we have  \eqref{autre condition sur Psi} if, for $x$ large enough:
\begin{displaymath}
\Psi (x) = \exp\big( \log x \,\log \log \log x\big).
\end{displaymath}
In fact, one has:
\begin{displaymath}
\frac{\Psi (Ax)}{\Psi (x)} =
\exp \bigg[ (\log x)\, \bigg( \log \frac{\log (\log x + \log A)}{\log \log x} \bigg) \bigg]  \,
\exp \big( \log A\, \log\log\log (Ax) \big),
\end{displaymath}
and $\exp \big( \log A\, \log\log\log (Ax) \big) = \big( \log \log (Ax) \big)^{\log A}$, whereas 
\begin{displaymath}
\log (\log x +\log A) 
= \log \log x + \log \Big( 1 + \frac{\log A}{\log x} \Big) 
\leq \log \log x + \frac{\log A}{\log x} \, \raise 1,5 pt \hbox{,}
\end{displaymath}
so:
\begin{displaymath}
\log \frac{\log (\log x + \log A)}{\log \log x} 
\leq \log \Big( 1 + \frac{\log A}{\log x \log \log x} \Big) 
\leq \frac{\log A}{\log x \log \log x} \,\cdot
\end{displaymath}
Hence:
\begin{displaymath}
\frac{\Psi (Ax)}{\Psi (x)} \leq \exp \Big( \frac{\log A}{\log \log x} \Big)\,\big( \log \log (Ax) \big)^{\log A} 
\lesssim \big( \log \log (Ax) \big)^{\log A} ,
\end{displaymath}
which gives \eqref{autre condition sur Psi} since, for every $\eps > 0$:
\begin{displaymath}
\big( \log \log (Ax) \big)^{\log A} = o\,\big( (\log x)^\eps \big) = o\,\big( [\log \Psi (x) ]^\eps \big).
\end{displaymath}
\medskip

On the other hand, if we require that the Orlicz function $\Psi$ grows sufficiently, we have:

\begin{proposition}\label{croissance Psi}
Assume that the Orlicz function $\Psi$ satisfies the condition:
\begin{equation}\label{condition croissance Psi}
\limsup_{x \to \infty} \frac{\Psi (Ax)}{[\Psi (x)]^2} > 0,
\end{equation}
for some $A > 1$.\par
Then, for every $p > 2$, there exists an analytic self-map $\phi \colon \D \to \D$ such that 
$C_\phi \colon H^2 \to H^2$ is in the Schatten class $S_p$, but $C_\phi \colon H^\Psi \to H^\Psi$ is not 
compact.
\end{proposition}

Note that condition~\eqref{condition croissance Psi} if not satisfied when $\Psi \in \Delta_2$, or for the 
Orlicz function $\Psi (x) = \e^{(\log (x+1) )^2} - 1$, \emph{i.e.} 
$\Psi (x) \approx \exp \big( (\log x)^2 \big)$.\par
\medskip

\noindent{\bf Proof.} Condition~\eqref{condition croissance Psi} implies that there exist some $\delta > 0$ and 
a sequence of positive numbers $x_n$ increasing to infinity such that:
\begin{displaymath}
\qquad \qquad \frac{\Psi (A x_n)}{[\Psi (x_n)]^2} \geq \delta, \qquad \text{for all } n \geq 1.
\end{displaymath}
If one sets $h_n = \frac{1}{\Psi (x_n)}\,$\raise 1,5pt \hbox{,} one has:
\begin{equation}\label{majoration par h carre}
h_n^2 \geq \frac{\delta}{\Psi \big[ A \Psi^{-1}( 1/h_n)\big]}\,\cdot
\end{equation}
Now, take $\alpha$ such that $\frac{2}{p} + 1 < \alpha < 2$, and consider the symbol $\phi =\phi_2$ 
constructed in \cite{JFA}, Theorem~4.1. Its composition operator 
$C_\phi \colon H^2 \to H^2$ belongs to $S_p$ (\cite{JFA},  Theorem~4.2). On the 
other hand, $C_\phi \colon H^\Psi \to H^\Psi$ is not compact, by \cite{CompOrli}, Theorem~4.18, since 
$\rho_\phi (h) \approx h^\alpha$ and hence $\rho_\phi (h_n)$ is not 
$\displaystyle o\,\bigg(\frac{1}{\Psi \big[ A \Psi^{-1}( 1/h_n)\big]}\bigg)$, 
by \eqref{majoration par h carre}. \hfill \vbox{}\qed
\bigskip

Condition~\eqref{condition croissance Psi} is a weaker growth condition than condition $\Delta^2$. Recall 
(see \cite{CompOrli}, \S~2.1) that the Orlicz function $\Psi$ satisfies the condition $\Delta^2$ if, for some 
$A > 1$, one has $[\Psi (x)]^2 \leq \Psi (Ax)$ for $x$ large enough. But condition $\Delta^2$ can be re-stated 
as:
\begin{equation}\label{Delta super 2 bis}
\liminf_{x \to \infty} \frac{\Psi (Ax)}{[\Psi (x)]^2} > 0.
\end{equation}
\par

Indeed, if  one has $[\Psi (x)]^2 \leq \Psi (Ax)$ for $x$ large enough, one obviously has 
$\liminf_{x \to \infty} \frac{\Psi (Ax)}{[\Psi (x)]^2} \geq 1$. Conversely, let 
$\delta = \liminf_{x \to \infty} \frac{\Psi (Ax)}{[\Psi (x)]^2}\,$. If $\delta > 1$, one clearly have 
$[\Psi (x)]^2 \leq \Psi (Ax)$ for $x$ large enough. If $0 < \delta \leq 1$, then, for $x$ large enough:
\begin{displaymath}
[\Psi (x)]^2 \leq \frac{2}{\delta} \Psi (Ax) \leq \Psi \Big(\frac{2}{\delta} Ax\Big)\,,
\end{displaymath}
by the convexity of $\Psi$, since $2/\delta \geq 1$. Hence $\Psi \in \Delta^2$ (with constant $2A /\delta$).
\par


\section{Orlicz functions with condition $\Delta^2$}\label{Delta super deux}

When $\Psi \in \Delta^2$, we have proved in \cite{CompOrli}, Corollary~3.26 a better result than the one in 
Proposition~\ref{croissance Psi}: \emph{Whenever $\Psi \in \Delta^2$, there exists a symbol 
$\phi \colon \D \to \D$ such that the composition operator $C_\phi \colon H^2 \to H^2$ is Hilbert-Schmidt, 
whereas $C_\phi \colon H^\Psi \to H^\Psi$ is not compact}.
\medskip

We have this improvement:

\begin{proposition}\label{proposition Schatten et Delta super deux} 
Let $\Psi$ be an Orlicz function satisfying condition $\Delta^2$. Then, there exists a symbol $\phi$ such 
that $C_\phi \colon H^2 \to H^2$ belongs to $S_p$ for every $p > 0$, but such that 
$C_\phi \colon H^\Psi \to H^\Psi$ is not compact.
\end{proposition}

\noindent{\bf Proof.} We shall use a \emph{lens map}. One can find a definition of such maps in 
\cite{Shap-livre}, page~27, but we shall use a slight variant of it: we shall use the same construction as in 
\cite{JFA}, Theorem~5.1 and Theorem~5.6, but with $f (z) = z^{1/2}$ instead of 
$f_\theta (z) = z (- \log z)^\theta$ or $f (z) = z \log ( - \log z)$. We shall recall this construction at the 
beginning of the next section. Note that 
$\Re (r\,\e^{i\alpha})^{1/2} = \sqrt{r}  \,\cos (\alpha /2)$ is positive ($| \alpha | < \pi /2$), so 
the associate symbol $\phi$ maps $\D$ into itself. \par
As in the proof of \cite{JFA}, Theorem~5.1, we need only to consider $f (it)$, with $t > 0$. But 
$f (it ) = \exp \big(\frac{1}{2} \log t + i \frac{\pi}{4} \big)$, so:
\begin{displaymath}
\Re f (it) = \sqrt{t}\cos (\pi /4) = \sqrt{t/2} \qquad \text{and} \qquad 
\Im f (it) = \sqrt{t}\sin (\pi/4) = \sqrt{t/2}.
\end{displaymath}
\par

Let us estimate $\rho_\phi (h)$. The condition $|\Re f (it)| \leq h$ implies $t \in [ -2h^2, 2h^2]$, and, 
just using the modulus constraint, we get $\rho_\phi (h) \lesssim h^2$. For the converse, if 
$t \in [ -2h^2, 2h^2]$, we have $|\Re f (it)| \leq h$ and $|\Im f (it) | \leq h$; then 
$\exp \big( f (it) \big) $ belongs to a window centered at $1$ and with size $Ch$. Hence 
$h^2 \lesssim \rho_\phi (h)$, and therefore $\rho_\phi (h) \approx h^2$.\par
It follows that $C_\phi \colon H^\Psi \to H^\Psi$ is not compact if $\Psi \in \Delta^2$. In fact, since 
$\Psi \in \Delta^2$, there is some $A_0 > 1$ such that $\Psi (A_0 x) \geq [\Psi (x)]^2$ for 
$x$ large enough. We get, with $x = \Psi^{-1} (1/h)$:
\begin{displaymath}
\frac{1}{\Psi \big(A_0 \Psi^{-1} (1/h) \big)}  \leq h^2.
\end{displaymath}
If $C_\phi$ were compact on $H^\Psi$, we should have (\cite{CompOrli}, Theorem~4.11):
\begin{displaymath}
\rho_\phi (h) = o \,\bigg( \frac{1}{\Psi \big(A_0 \Psi^{-1} (1/h) \big)} \bigg)\,,
\end{displaymath} 
which is not the case.\par
\smallskip

On the other hand, $C_\phi$ is in $S_p$ for every $p > 0$. Indeed, $\phi (\e^{it})$ is in the 
dyadic Carleson window: 
\begin{displaymath}
W_{n,j}=\Big\{ z\in\D \,;\  1 - 2^{-n} \leq |z| < 1\, ,\quad 
\frac{2j\pi}{2^n} \leq \arg(z) < \frac{2 (j + 1) \pi}{2^n} \Big\},
\end{displaymath}
($j=0, 1, \ldots, 2^n - 1$, $n=1,2,\ldots$) if and only if, for $h_n =2^{-n}$, one has, up to constants:
\begin{displaymath}
t \leq 2h_n^2 \qquad \text{and} \qquad j^2 h_n^2 \leq t \leq (j + 1)^2 h_n^2. 
\end{displaymath}
This is possible for at most a fixed number, say $N$, of values of $j$; hence:
\begin{align*}
\sum_{n=1}^\infty \sum_{j=0}^{2^n -1} 2^{n p/2} m_\phi ( W_{n,j} )^{p/2} 
& \leq N \sum_{n=1}^\infty 2^{n p/2} \big( \rho_\phi (2^{-n}) \big)^{p/2} \\
& \leq N \sum_{n=1}^\infty 2^{n p/2} 2^{-2n p/2} = N \sum_{n=1}^\infty 2^{- n p/2} < +\infty.
\end{align*}
Then Luecking's theorem (\cite{Lue}) and \cite{JFA}, Proposition~3.3 imply that $C_\phi \in S_p$. \qed
\bigskip

Now, in the opposite direction, one has proved in \cite{CompOrli}, Theorem~3.24, that if $C_\phi$ is 
compact on $H^\Psi$, with $\Psi \in \Delta^2$, then $C_\phi \in S_p (H^2)$ for every $p > 0$ (though we 
only stated that $C_\phi \in S_1$). We have the following improvement.

\begin{proposition}\label{Delta super deux et Carleson}
Let $\Psi$ be an Orlicz function satisfying condition $\Delta^2$. Assume that $C_\phi$ is compact on 
$H^\Psi$. Then $m_\phi$ is an $\alpha$-Carleson measure, for every $\alpha \geq 1$, and hence 
$C_\phi \in S_p$ for every $p > 0$.
\end{proposition}

\noindent{\bf Proof.} Since $\Psi \in \Delta^2$, there is some $A_0 > 1$ such that 
$\Psi (A_0 x) \geq [\Psi (x)]^2$, for $x$ large enough. Hence, for every positive integer $n$, one has 
$\Psi (A_0^n x) \geq [\Psi (x)]^{2^n}$. Taking $x = \Psi^{-1} (1/h)$\,, one gets:
\begin{displaymath}
\Psi \big( A_0^n \Psi^{-1} (1/h) \big) \geq \frac{1}{h^{2^n}}\,\cdot
\end{displaymath}
Now, if $C_\phi$ is compact on $H^\Psi$, one has, by \cite{CompOrli}, Theorem~4.11,~1), with $A = A_0^n$:
\begin{displaymath}
\rho_\phi (h) \leq \frac{1}{\Psi \big( A_0^n \Psi^{-1} (1/h) \big)}\,;
\end{displaymath}
therefore $\rho_\phi (h) \leq h^{2^n}$, and $m_\phi$ is a $2^n$-Carleson measure.\par
The last assertion now follows from \cite{JFA}, Corollary~3.2. \qed


\section{Case of the symbols $\phi_\theta$} \label{phi indice theta dans H Psi}

In \cite{JFA}, Section~5, we considered composition operators $C_{\phi_\theta}$ which are, for 
every $\theta > 0$, compact on $H^2$, and hence on $H^p$ for every $p < \infty$. We are now going to 
examine when these composition operators are compact on $H^\Psi$.\par
\smallskip

Let us recall how the $\phi_\theta$'s are constructed.\par
For $\Re z > 0$, consider the principal determination of the logarithm $\log z$. For $\eps >0$ small enough, 
the function $f_\theta (z) = z (- \log z )^\theta$ has strictly positive real part on 
$V_\eps = \{ z\in \C\,;\ \Re z >0 \ \text{and}\ |z | < \eps\}$ and 
$\Re f_\theta^\ast ( z) > 0$ for all $z \in \partial V_\eps \setminus\{0\}$. Now, if 
$g_\theta$ is the conformal mapping from $\D$ onto $V_\eps$, which maps $\T = \partial \D$ 
onto $\partial V_\eps$, and with $g_\theta (1) = 0$ and $g'_\theta (1) = -\eps/4$, we set:
\begin{displaymath}
\phi_\theta = \exp ( - f_\theta \circ g_\theta).
\end{displaymath}
$\phi_\theta$ maps $\D$ into itself and $|\phi_\theta^\ast| < 1$ on $\partial \D \setminus\{1\}$, and we 
proved (\cite{JFA}, Proposition~5.3) that:
\begin{equation}\label{rho pour phi indice theta}
\rho_{\phi_\theta} (h) \approx \frac{h}{(\log 1/h )^\theta}\cdot
\end{equation}
\bigskip

Now, recall (\cite{CompOrli}, Theorem~4.18) that $C_{\phi_\theta}$ is compact on $H^\Psi$ if and 
only if, for every $A > 0$: 
\begin{displaymath}
\rho_{\phi_\theta} (h) = o\bigg( \frac{1}{\Psi [ A \Psi^{-1} (1/h)]}\bigg)\,, \qquad h \to 0.
\end{displaymath}
By \eqref{rho pour phi indice theta}, that gives, with $u = \Psi^{-1}(1/h)$:

\begin{theoreme}\label{Psi delta 1}
The composition operator $C_{\phi_\theta}$ is compact on $H^\Psi$ if and only if for every $A > 1$, one has:
\begin{equation}\label{compacite phi indice theta sur H Psi}
\qquad \qquad \Psi (Au) = o\,\big( \Psi (u) [\log \Psi (u)]^\theta\big)\,,\qquad \text{as}\ u \to \infty.
\end{equation}
\end{theoreme}

For example, for every $\theta > 0$, $C_{\phi_\theta}$ is compact  on $H^\Psi$ if 
$\Psi (x) \approx x^{\log \log \log x} = \e^{(\log x)(\log \log \log x)}\!$, but is not compact on 
$H^\Psi$ if $\Psi (x) \approx x^{\log \log x} = \e^{(\log x)(\log \log x)}$.
\par\smallskip

Remark that this condition~\eqref{compacite phi indice theta sur H Psi} is almost the same as 
condition~\eqref{autre condition sur Psi}: in \eqref{compacite phi indice theta sur H Psi} one is requiring the 
condition to be satisfied for some $\theta > 0$, whereas in \eqref{autre condition sur Psi} one is 
requiring it to be satisfied for every $\eps > 0$ (see however the remark at the end of this Section). The 
\emph{little-oh} assumption in \eqref{compacite phi indice theta sur H Psi} is actually not very important: 
suppose that, for every $ A > 1$, one has $\Psi (Au)  \leq C\,\Psi (u) [\log \Psi (u)]^\theta$ for $u \geq u_A$ 
(we need a uniform constant $C > 0$), then, for every $ A > 1$, the convexity of $\Psi$ gives, for every 
$\eps > 0$,  $\Psi (A u) \leq \Psi (A_\eps u)$, with $A_\eps = A/\eps$; hence we get 
$\Psi (Au) \leq C\,\eps\, \Psi (u) [\log \Psi (u)]^\theta$, and the \emph{little-oh} property follows.
\par
\medskip

Note that, since $C_{\phi_\theta}$ is in $S_p$ for $p > 4/\theta$ (\cite{JFA}, Proposition~5.3),  
Theorem~\ref{Psi delta 1} is an improvement of Proposition~\ref{croissance Psi}, since 
\eqref{compacite phi indice theta sur H Psi} is not satisfied if and only if there exists some $A >1 $ such that 
$\limsup_{x \to \infty} \frac{\Psi (Ax)}{\Psi (x) [\log \Psi (x)]^\theta} > 0$, which is of course implied 
by condition \eqref{condition croissance Psi}:

\begin{theoreme}\label{meilleure condition sur Psi}
Let $p > 0$ and $\theta > 4/ p$. Assume that the Orlicz function $\Psi$ satisfies, for some $A > 1$, the 
condition:
\begin{equation}\label{Delta 1, theta}
\limsup_{x \to \infty} \frac{\Psi (Ax)}{\Psi (x) [\log \Psi (x)]^\theta} > 0. 
\end{equation} 
Then the composition operator $C_{\phi_\theta}$, which is in $S_p (H^2)$, is not compact on $H^\Psi$.
\par
If the Orlicz function $\Psi$ satisfies condition $\Delta^1$, then condition~\eqref{Delta 1, theta} is 
satisfied. 
\end{theoreme}

Remark that, for $\theta > 2$, Theorem~\ref{meilleure condition sur Psi} improves  
Corollary~3.26 of \cite{CompOrli}. Recall that the Orlicz function $\Psi$ satisfies the condition 
$\Delta^1$ if there is a constant $A > 1$ such that $x \Psi (x) \leq \Psi (A x)$ for $x$ large enough 
(see \cite{CompOrli}, \S~2.1). \par
\medskip

\noindent{\bf Proof of Theorem~\ref{meilleure condition sur Psi}.} The first part was explained before 
the statement. For the second part, assume that, for every $A > 1$,  
condition~\eqref{compacite phi indice theta sur H Psi} is satisfied; in particular, we should have, for some 
$u_0 > 0$:
\begin{equation}\label{inegalite Psi}
\qquad \qquad \Psi (Au) \leq \Psi (u) \big( \log \Psi (u) \big)^\theta, \qquad u \geq u_0.
\end{equation}
Let
\begin{displaymath}
p_n = \Psi (A^n) \qquad \text{and} \qquad l_n = \log p_n.
\end{displaymath}
One has then, by \eqref{inegalite Psi}, for $n$ large enough:
\begin{displaymath}
p_{n+1} \leq p_n (\log p_n)^\theta \qquad \text{and} \qquad l_{n+1} \leq l_n + \theta \log l_n.
\end{displaymath}
Let
\begin{displaymath}
h (t) =\int_2^t \frac{dx}{\theta \log x}\cdot
\end{displaymath}
One has, since $h'$ decreases:
\begin{displaymath}
0  < h (l_{n+1}) - h (l_n) \leq (l_{n+1} - l_n) h' (l_n) \leq (\theta \log l_n) \frac{1}{\theta \log l_n} 
\leq 1.
\end{displaymath}
Hence $h (l_n) \lesssim n$, \emph{i.e.}:
\begin{displaymath}
\int_2^{l_n} \frac{dx}{\theta \log x} \lesssim n.
\end{displaymath}
Since this integral is $\gtrsim \frac{l_n}{\theta \log l_n}$, we get:
\begin{equation}\label{inegalite log Psi}
\frac{\log \Psi (A^n)}{\theta \log \log \Psi (A^n)} \lesssim n.
\end{equation}
Hence:
\begin{displaymath}
\log \log \Psi (A^n) - \log \log \log \Psi (A^n) \lesssim \log n,
\end{displaymath}
and so:
\begin{displaymath}
\limsup_{n \to \infty} \frac{\log \log \Psi (A^n)}{\log n} \leq 1.
\end{displaymath}
\par

On the other hand, $\Psi (A^n) \geq A^n$, \emph{i.e.} $\log \Psi (A^n) \geq n \log A$. \par
Therefore $\log \log \Psi (A^n) \sim \log n$, and hence, using \eqref{inegalite log Psi}:
\begin{displaymath}
\log \Psi (A^n) \lesssim  2\theta n \log n.
\end{displaymath}
But this prevents $\Psi$ from satisfying the condition $\Delta^1$: if $\Psi \in \Delta^1$, there exists 
$A > 1$ such that $u \Psi (u) \leq \Psi (A u)$ , for $u$ large enough, and this implies, for $n$ large 
enough, that $A^n \Psi (A^n) \leq \Psi (A^{n+1})$, and hence 
$A^{n(n+1)/2} \Psi (A) \lesssim \Psi (A^n)$, so that $\log \Psi (A^n) \gtrsim n^2$. \qed
\bigskip

\noindent{\bf Remark.} As a consequence of Theorem~\ref{Psi delta 1}, one has the following striking result. 

\begin{corollaire}\label{pas d'interpolation}
There exists an Orlicz function $\Psi$ such that $H^{3+\eps} \subseteq H^\Psi \subseteq H^3$ for every 
$\eps >0$, and a composition operator $C_\phi$ which is compact on $H^3$ and on $H^{3+\eps}$, but 
not compact on $H^\Psi$.
\end{corollaire}

\noindent{\bf Remark.} It is known (\cite{Shap-Tay}, Theorem 6.1), that the compactness of $C_\phi$ on 
$H^{p_0}$ for some $p_0 <\infty$ implies its compactness on $H^p$ for all $p < \infty$; this result follows 
from the Riesz factorization Theorem. It follows from Corollary~\ref{pas d'interpolation} that there is no 
such factorization for Hardy-Orlicz spaces in general.
\medskip

\noindent{\bf Proof.} Consider the symbol $\phi = \phi_\theta$. 
The compactness of  $C_\phi$ from $H^2$ into itself implies its compactness from 
$H^p$ into itself for every $p \geq 1$; in particular its compactness on $H^3$ and on $H^{3+\eps}$.\par 
Let now $\Psi$ the Orlicz function contructed in \cite{critere}. This function is such that 
$x^3/3 \leq \Psi (x)$ for all $x\geq 0$, but $\Psi (n!) \leq (n!)^3$ for all positive integer $n \geq 1$. 
Condition~\eqref{compacite phi indice theta sur H Psi} is not satisfied when $\theta < 1$ because 
$\Psi (3.k!) \geq k.(k!)^3$ and 
\begin{displaymath}
\Psi (k!) [\log \Psi (k!)]^\theta \leq (k!)^3[ 3 \log (k!)]^\theta \lesssim (k!)^3 (k \log k)^\theta;
\end{displaymath}
hence $C_{\phi_\theta}$ is not compact on $H^\Psi$ if $\theta < 1$.\qed
\medskip

\noindent{\bf Remark.} Condition \eqref{compacite phi indice theta sur H Psi} depends on $\theta$, but 
if $\Psi$ is \emph{regular}, the compactness of $C_{\phi_\theta}$ on $H^\Psi$ does not actually depend on 
$\theta$, because, writing:
\begin{displaymath}
[\log \Psi (u)]^\theta \gtrsim \frac{\Psi (Au)}{\Psi (u)} 
= \frac{\Psi (Au)}{\Psi (\sqrt{A} u)} \frac{\Psi (\sqrt{A} u)}{\Psi (u)} \,\raise 1,5pt \hbox{,}
\end{displaymath}
one gets:
\begin{displaymath}
\frac{\Psi (\sqrt{A} u)}{\Psi (u)} \lesssim [\log \Psi (u)]^{\theta/2} 
\qquad \text{or else} \qquad 
 \frac{\Psi (A u)}{\Psi (\sqrt{A} u)} \lesssim [\log \Psi (u)]^{\theta/2} \cdot
\end{displaymath}
Hence, if $\Psi$ satisfies the following condition of regularity:
\begin{equation}\label{nabla}
\frac{\Psi ( \sqrt{A} u)}{\Psi (u)} \lesssim \frac{\Psi (Au)}{\Psi ( \sqrt{A} u)}\,\raise 1,5pt 
\hbox{,} \qquad \forall A > 1, \forall u > 0,
\end{equation}
we get, for every $A >1$ and every $u >0$:
\begin{displaymath}
\frac{\Psi (\sqrt{A} u)}{\Psi (u)} \lesssim [\log \Psi (u)]^{\theta/2} \,\raise 1,5pt \hbox{,}
\end{displaymath}
and condition \eqref{compacite phi indice theta sur H Psi} remains true by replacing 
$\theta$ with $\theta/2$. \par
Note that condition~\eqref{nabla} is satisfied whenever the Orlicz function $\Psi$ satisfies the condition 
$\nabla_0$ defined in \cite{CompOrli}, Definition~4.5; indeed, by Proposition~4.6 of \cite{CompOrli}, $\Psi$ 
satisfies $\nabla_0$ if and only if  there exists some $u_0 >0$ such that, for every $\beta > 1$, there exists 
$C_\beta \geq 1$ such that:
\begin{displaymath}
\frac{\Psi (\beta u)}{\Psi (u)} \leq \frac{\Psi (\beta C_\beta v)}{\Psi (v)}\,\raise 1,5 pt \hbox{,} 
\qquad u_0 \leq u \leq v.
\end{displaymath}
This condition implies \eqref{nabla}, by taking  $\beta = \sqrt A$ and $v = \sqrt{A}u$, since, by convexity, 
$\Psi (\beta C_\beta v) \leq C_\beta \Psi (\beta v)$.\par
When we impose in condition~\eqref{nabla} an inequality with factor $1$, for $u$ large enough, one has:
\begin{displaymath}
\big[\Psi ( \sqrt{A} u)\big]^2 \leq \Psi (u)\, \Psi (Au)\,, \qquad \forall A > 1, \forall u \geq u_0,
\end{displaymath}
and if one sets $u = \e^x$ and $Au = \e^y$, we get:
\begin{displaymath}
\log \Psi \Big[ \exp \Big( \frac{x + y}{2} \Big)\Big] 
\leq \frac{1}{2} \big[\log \Psi (\e^x) + \log \Psi (\e^y)\big]\,,\qquad x_0 \leq x < y,
\end{displaymath}
which means that the function $\kappa (x) = \log \Psi (\e^x)$ is convex for $x$ large enough. By 
\cite{CompOrli}, Proposition~4.7, that means that the Orlicz function $\Psi$ satisfies the condition 
$\nabla_0$ with constant $1$. This regularity condition is satisfied, for example if 
$\Psi (x) = x^p$, $p\geq 1$, or $\Psi (x) = \e^{(\log (x+1))^\alpha} - 1$, or else 
$\Psi (x) = \e^{x^\alpha} -1$, $\alpha \geq 1$.


\section{Comparison with $H^\infty$}\label{Other results}

It is well-known, and easy to see, that $C_\phi$ is compact on $H^\infty$ if and only if $\| \phi\|_\infty < 1$. 
For composition operators, $H^\infty$ is a very special case; indeed, one has:

\begin{proposition}
\hfill\par
1) Given any analytic self-map $\phi \colon \D \to \D$ such that $\| \phi\|_\infty = 1$, there exists an Orlicz 
function $\Psi$ such that $C_\phi$ is not compact on $H^\Psi$.\par
2) Given any Orlicz function $\Psi$, there exists an analytic self-map $\phi \colon \D \to \D$ such that 
$\| \phi\|_\infty = 1$ and for which $C_\phi$ is compact on $H^\Psi$.
\end{proposition}

\noindent{\bf Proof.} 1) Since $\| \phi\|_\infty = 1$, one has $\rho_\phi (h) > 0$ for every $0 < h < 1$.\par
We construct piecewise a strictly convex function $\Psi \colon [0,+\infty[ \to \R_+$ such that:
\begin{displaymath}
\Psi (2^{n+1}) \geq \frac{1}{\rho_\phi \big( 1 / \Psi (2^n)\big)} \, \raise 1,5 pt \hbox{,}
\end{displaymath}
Then, if $h_n = 1/ \Psi (2^n)$, one has $h_n \to 0$, and
\begin{displaymath}
\rho_\phi (h_n) \geq \frac{1}{\Psi \big( 2 \Psi^{-1} (1/h_n) \big) }\cdot
\end{displaymath}
It follows from \cite{CompOrli}, Theorem~4.18, that $C_\phi$ is not compact on $H^\Psi$.\par
2) Let $a \in M^\Psi (\T)$, the Morse-Transue space (see \cite{CompOrli}) a positive function such that 
$a \notin L^\infty (\T)$ and such that $\Psi \circ a \geq 2$, and let:
\begin{displaymath}
h = 1 -\frac{1}{\Psi \circ a}\cdot
\end{displaymath}
Since $a$ is not essentially bounded, one has $\| h\|_\infty =1$. Moreover, one has $h \geq 1/2$, so 
$\log h$ is integrable; we can hence define the outer function:
\begin{displaymath}
\phi (z) =\exp\bigg[ \int_\T \frac{ u + z}{u - z} \log h (u)\,dm (u) \bigg]\,, \qquad | z | < 1.
\end{displaymath}
One has $| \phi^\ast (\e^{i\theta})| = h (\e^{i\theta})$, so $\phi$ is an analytic self-map from $\D$ into 
$\D$, and $\| \phi \|_\infty = 1$. Since, for every $A >0$:
\begin{displaymath}
\Psi \bigg[A \Psi^{-1} \Big( \frac{1}{1 -|\phi^\ast|} \Big)\bigg] = \Psi (A a) 
\end{displaymath}
is integrable, because $a \in M^\Psi (\T)$, the composition operator $C_\phi$ is $M^\Psi$-order bounded 
(see \cite{CompOrli}, Proposition~3.14), and hence is compact on $H^\Psi$. \qed
\bigskip


\bigskip

\vbox{\noindent{\it 
{\rm Pascal Lef\`evre}, Universit\'e d'Artois,\\
Laboratoire de Math\'ematiques de Lens EA~2462, \\
F\'ed\'eration CNRS Nord-Pas-de-Calais FR~2956, \\
Facult\'e des Sciences Jean Perrin,\\
Rue Jean Souvraz, S.P.\kern 1mm 18,\\ 
62\kern 1mm 307 LENS Cedex,
FRANCE \\ 
pascal.lefevre@euler.univ-artois.fr 
\smallskip

\noindent
{\rm Daniel Li}, Universit\'e d'Artois,\\
Laboratoire de Math\'ematiques de Lens EA~2462, \\
F\'ed\'eration CNRS Nord-Pas-de-Calais FR~2956, \\
Facult\'e des Sciences Jean Perrin,\\
Rue Jean Souvraz, S.P.\kern 1mm 18,\\ 
62\kern 1mm 307 LENS Cedex,
FRANCE \\ 
daniel.li@euler.univ-artois.fr
\smallskip

\noindent
{\rm Herv\'e Queff\'elec},
Universit\'e des Sciences et Technologies de Lille, \\
Labo\-ratoire Paul Painlev\'e U.M.R. CNRS 8524, \\
U.F.R. de Math\'ematiques,\\
59\kern 1mm 655 VILLENEUVE D'ASCQ Cedex, 
FRANCE \\ 
queff@math.univ-lille1.fr
\smallskip

\noindent
{\rm Luis Rodr{\'\i}guez-Piazza}, Universidad de Sevilla, \\
Facultad de Matem\'aticas, Departamento de An\'alisis Matem\'atico,\\ 
Apartado de Correos 1160,\\
41\kern 1mm 080 SEVILLA, SPAIN \\ 
piazza@us.es\par}
}

\end{document}